\documentclass[12pt,a4paper]{article}
  
\usepackage[english]{babel}
\usepackage{amsmath,amssymb} 
\usepackage{graphicx} 
\usepackage{cite} 
\usepackage{url} 
  
\def\mbf#1{\protect{\makebox{\boldmath $#1$}}}

\newcommand{\mbR}{\mathbb{R}} 
\newcommand{\mboR}{\overline{\mbR}} 
\newcommand{\mbN}{\mathbb{N}} 
\newcommand{\mbIR}{\mathbb{IR}} 
\newcommand{\mbIoR}{\mathbb{I}\overline{\mbR}} 
\newcommand{\mbKR}{\mathbb{KR}} 
\newcommand{\mbKoR}{\mathbb{K}\overline{\mbR}} 
  
\DeclareMathOperator{\Mid}{mid} 
\DeclareMathOperator{\rad}{rad} 
\DeclareMathOperator{\dual}{dual} 
\DeclareMathOperator{\pro}{pro} 
\DeclareMathOperator{\ext}{extr} 
\DeclareMathOperator{\abs}{abs} 
\newlength{\lvee} 
\newcommand{\W}{\settowidth{\lvee}{$\bigvee$}\bigvee\hspace{-\lvee}\bigwedge\nolimits} 
  
\newcommand{\QAb}{Q(\bA,\bb,\calA,\beta)}
\newcommand{\Qi}{Q_{i:}(\bA,\bb,\calA,\beta)}
\newcommand{\QAEAb}{Q^{\forall\exists}(\bA,\bb,\calA,\beta)}
\newcommand{\Qsigma}{Q^{\sigma}(\bA,\bb,\calA,\beta)} 
\newcommand{\Ac}{\mbf{A}^{\mathfrak{c}}}
\newcommand{\bc}{\mbf{b}^{\mathfrak{c}}}
  
\newcommand{\un}{\underline}
\newcommand{\ov}{\overline}
  
\newcommand{\bA}{\mbf{A}}
\newcommand{\calA}{{\cal A}}
\newcommand{\cAs}{{\cal A}^s}
\newcommand{\betas}{\beta^s}
\newcommand{\cA}{\check{\mbf{A}}}
\newcommand{\hA}{\hat{\mbf{A}}}
  
\newcommand{\bB}{\mbf{B}}

\newcommand{\bb}{\mbf{b}} 
\newcommand{\cb}{\check{\mbf{b}}} 
\newcommand{\hb}{\hat{\mbf{b}}}

\newcommand{\bu}{\mbf{u}}
\newcommand{\cu}{\check{\mbf{u}}}
\newcommand{\hu}{\hat{\mbf{u}}}
  
\newcommand{\bv}{\mbf{v}}
\newcommand{\cv}{\check{\mbf{v}}}
\newcommand{\hv}{\hat{\mbf{v}}}
  
\newcommand{\bw}{\mbf{w}}
  
\newcommand{\conj}[1]{\underset{#1}{\scalebox{1.5}{\&}}\ } 
  
\newtheorem{theorem}{Theorem} 
\newtheorem{corollary}{Corollary} 
  
\begin{title}
{\bf Quantifier-free descriptions \\[3pt] 
     for quantifier solutions     \\[3pt] 
     to interval linear systems of relations\footnote{The work is a revised and extended  
     English translation of the original Russian work \textsc{Sharaya, I.A.} Quantifier-free
     descriptions of interval-quantifier linear systems // \textsl{Trudy Instituta Matematiki 
     i Mekhaniki UrO RAN (Proceedings of the Institute of Mathematics and Mechanics, Ural 
     Branch of the Russian Academy of Sciences)}. -- 2014. -- Vol.~20, No.~2. -- P.~311--323.} 
     \\[3mm]} 
\end{title} 
  
\author{\Large\sf Irene A. Sharaya} 
\date{\normalsize Institute of Computational Technologies SB RAS \\[2pt] 
      and Novosibirsk State University, \\[2pt] 
      Novosibirsk, Russia} 
  
\begin{document}
\maketitle 
  
\begin{abstract} 
We study systems of relations of the form $Ax\,\sigma\,b$, where $\sigma$ is a vector 
of binary relations with the components ``$=$'', ``$\geq$'' and ``$\leq$'', and the 
parameters (elements of the matrix $A$ and right-hand side vector $b$) can take values 
from prescribed intervals. What is considered to be the set of its solutions depends 
on which logical quantifier is associated with each interval-valued parameter and 
what is the order of the quantifier prefixes for certain parameters. For solution sets 
that correspond to the quantifier prefix of a general form, we present equivalent 
quantifier-free descriptions in the classical interval arithmetic, in Kaucher 
comp\-lete interval arithmetic and in the usual real arithmetic. \\[2mm] 
\end{abstract} 
  
  
\section{Introduction} 
  
  
\subsection{Quantifier solutions to interval linear systems}

In the classical interval arithmetic $\mbIR$, an interval is a non-empty bounded 
connected closed subset of the real line $\mbR$. According to the notation standard 
\cite{INotation}, we will denote interval objects in bold type ($\mbf{A}$, $\mbf{B}$, 
\ldots, $\mbf{y}$, $\mbf{z}$), in contrast to usual point (non-interval) quantities 
that are not specifically distinguished. 
  
We consider a system of linear equations and inequalities of the form 
\begin{equation*} 
Ax\,\sigma\,b, \qquad A\in\mbR^{m\times n},\ x\in\mbR^n,\ b\in\mbR^m,\
\sigma\in\{=,\geq,\leq\}^m,\ \ m,n\in\mbN,
\end{equation*} 
where $x$ is a vector of unknowns, $\sigma$ is a vector of binary relations, with 
the components ``$=$'', ``$\geq$'' and ``$\leq$'', and every parameter $u\in\mbR$ 
(which may be an element of the matrix $A$ or of the right-hand side $b$) can take 
values within the prescribed eponymous interval $\bu$ from $\mbIR$. 
  
As an example of such ``mixed'' systems, we can consider a $3\times 2$-system 
\begin{equation*} 
\arraycolsep=2pt 
\left\{ 
\begin{array}{rcl} 
2x_{1} - 3x_{2} &\leq& 4, \\[2pt]  
5x_{1} + 6x_{2} & =  & 7, \\[2pt]
-x_{1} + 4x_{2} &\geq& 5. 
\end{array} 
\right. 
\end{equation*} 
  
With each parameter $u$, we connect either the universal quantifier ``$\forall$'' 
or the existential quantifier ``$\exists$'' as well as the corresponding elementary 
quantifier prefix $(\forall u\in\bu)$ or $(\exists u\in\bu)$. Such interval uncertainty 
of parameters can be specified by the interval matrix $\bA\in\mbIR^{m\times n}$, 
the quantifier matrix $\calA$ of the same size as $\bA$, the interval vector 
$\bb\in\mbIR^m$ and the quantifier vector $\beta$ of length $m$. All the elementary 
quantifier prefixes can be written down in an arbitrary order, and we denote 
the resulting prefix of the length $m(n+1)$ as $\QAb$. 
  
\bigskip\noindent
\textbf{Definition 1}\, 
For given interval matrix $\mbf{A}$, interval vector $\mbf{b}$, and quantifier 
prefix $\QAb$, the \textit{interval-quantifier linear system of relations}, or 
\textit{interval-quantifier linear system} in short, will be called the predicate 
of the form $\QAb (Ax\,\sigma\,b)$. A vector $\tilde{x}\in\mbR^n$ is referred to as 
\textit{solution} to the interval-quantifier linear system if the predicate 
$\QAb (Ax\,\sigma\,b)$ takes the value ``true'' in $\tilde{x}$. 
  
\bigskip 
This definition is, in fact, a further refinement of the general ideas expressed 
in \cite{SharyRC02}. Also, the above introduced interval-quantifier linear systems 
are a natural generalization of interval linear systems, which have long been studied 
in interval analysis. \emph{Interval linear system} of the form $\bA x\,\sigma\,\bb$ 
is a formal record, for which we specially stipulate what is considered a solution 
in each specific case. Usually, interval linear systems of only equations or of only 
inequalities of the same sign are considered, while their solutions are taken 
as formal (algebraic) solutions, AE-solutions, strong solutions, weak solutions, 
tolerable solutions, controllable solutions (sometimes called simply ``control solutions''), 
etc. (see \ [\citen{SharyBook,SharyRC02}; \citen{FNRRZ}, chapter~2] and references 
in these publications). 
  
In order to agree with the existing terminology in this field, the solutions to the 
interval-quantifier linear system $\QAb (Ax\,\sigma\,b)$ will also be called the 
\emph{quantifier} solutions of the interval linear system $\bA x\,\sigma\,\bb$. 
Notice that AE-solutions, strong solutions, weak solutions, tolerable solutions, 
controllable solutions to the system $\bA x\,\sigma\,\bb$ are subsumed under 
the quantifier solutions. In particular, the AE-solutions to interval linear systems 
of equations are their quantifier solutions for which the determining predicate 
$\QAb (Ax\,\sigma\,b)$ has the so-called AE-form, that is, in which all occurrences 
of the universal quantifier ``$\forall$'' precede the occurrences of the existential 
quantifier ``$\exists$''. 
  
The notation $\QAb (Ax\,\sigma\,b)$ defines all possible interval-quantifier linear 
systems in parametric form. The parameters of the description are $\bA$, $\bb$, 
$\calA$, $\beta$, $\sigma$ and the order of elementary quantifier prefixes in $Q$ 
(since the elementary prefixes with different quantifiers do not always commute). 
Imposing additional constraints on the parameters, we obtain different classes (subsets) 
of interval-quantifier linear systems. For example, if we require equality to be the 
value of each component of the relation vector $\sigma$, then we obtain a class 
of interval-quantifier systems of linear equations.

  
\subsection{Transition to quantifier-free descriptions}

Interval-quantifier linear systems and their solutions were introduced in the previous 
section through a logical predicate of the first order. Predicative description is close 
to the formulation of practical problems, but it allows very limited means of theoretical 
investigation and is not at all suitable for calculations. As a consequence, 
the following problem arises: 
  
\bigskip\noindent
\textbf{Problem}\ \  
For the widest possible subset of interval-quantifier linear systems, find a convenient 
quantifier-free description of their solutions in algebraic systems with sufficiently 
developed tools for equivalent transformations, study and computation. 
  
\bigskip 
Usually, the solution sets to interval systems of equations and inequalities are described 
using real arithmetic in $\mbR$ [\citen{FNRRZ}, chapter 2;\, \citen{Eryomin88}, pp.~93--95; 
\, \citen{OettliPrager,Gerlach,Vatolin84,Rohn85,LakeyevNoskov94,Lakeyev03}], since it is 
simple, familiar, has good properties, and we can apply developed numerical methods 
in $\mbR$. For various subclasses of interval quantifier systems of linear equations, 
a number of quantifier-free descriptions have been obtained in interval arithmetics 
\cite{Beeck,Neumaier86,SharyBook,SharyRC02}. Despite poor algebraic properties of 
the interval arithmetics (the absence of the distributivity, etc.), these descriptions 
turned out to be very useful. For instance, the description of the AE-solution sets 
of interval systems of linear equations made it possible to construct a general theory 
of these solutions and interval numerical methods for inner and outer estimation of 
the AE-solution sets \cite{SharyBook,SharyRC02}. 
  
The features of the quantifier-free descriptions proposed in this paper are as follows: 
\begin{itemize} 
\item[1.] 
They expand the class of the interval-quantifier linear systems for which a convenient 
description of solutions can be given in comparison with those known descriptions where 
non-negativity of $x$ is not required.  (The nonnegativity requirement on the vector 
of unknowns can be formulated as an additional restriction on the parameters $\bA$, 
$\bb$ and $\sigma$. A.\,Vatolin in \cite{Vatolin84} obtained quantifier-free 
descriptions for solutions of general interval-quantifier linear systems, but his result 
is only valid under nonnegativity condition on the variables, which is quite severe 
limitation in practice. The class $Q^\sigma$ of interval-quantifier linear systems 
we discuss in the present paper has no restrictions on $\bA$, $\bb$ and $\sigma$, but 
it has a restriction on the order of the elementary quantifier prefixes.) 
\item[2.] 
Quantifier-free descriptions of the solutions are obtained in ordinary real arithmetic 
$\mbR$, classical interval arithmetic $\mbIR$, and in Kaucher interval arithmetic $\mbKR$. 
This enables us to carry out investigation of the solution sets and computation 
with them by both real and interval methods. 
\end{itemize}

  
\section{Necessary facts about interval arithmetics}

In this section, we give the necessary information on various interval arithmetics. 
The desire to improve the properties of the classical interval arithmetic $\mbIR$ led 
to the appearance of its various extensions. The most popular of them is \emph{Kaucher 
interval arithmetic} $\mbKR$ developed by E.\,Kaucher \cite{Kaucher80}. E.\,Garde\~{n}es 
and A.\,Trepat \cite{GardTr80} and then S.\,Markov \cite{Markov95} proposed another 
similar constructions. All these researchers constructed extensions of the classical 
interval arithmetic $\mbIR$ on the basis of different principles, which was reflected 
in the names of the corresponding algebraic structures: extended interval arithmetic 
\cite{Kaucher80}, modal interval arithmetic \cite{GardTr80}, arithmetic of directed 
intervals \cite{Markov95}. However, despite the difference in their construction, 
all three algebraic systems coincide up to notation. 
  
\textit{Interval} in $\mbKR$ is just a record of the form $[a,b]$, where $a,b\in\mbR$. 
In $\mbIR$, the values $a$ and $b$ should additionally satisfy the requirement $a\leq b$. 
Intervals are also denoted by small boldface letters, e.\,g., $\bu\in\mbKR$. If $\bu$ 
and $[a,b]$ denote the same interval, then $a$ is called left (lower) endpoint of the 
interval, which is written as $\un{\bu}$, and $b$ is called the right (upper) endpoint 
of the interval $\bu$, which is written as $\ov{\bu}$. Therefore, $\bu\equiv[\un{\bu}, 
\ov{\bu}]$. The intervals from $\mbIR$, as was stated in Introduction, can be considered 
as subsets of the real axis $\mbR$: 
\begin{equation*} 
[\un{\bu},\ov{\bu}] = \{\,u\in\mbR\mid \un{\bu}\leq u\leq\ov{\bu}\,\}. 
\end{equation*} 
  
In this paper, we will mainly use concepts and properties of Kaucher interval arithmetic, 
and we present them below. 
  
Two intervals are considered equal if both their left and right endpoints coincide: 
\begin{equation*} 
\bu = \bv \quad 
\stackrel{\mathrm{def}}{\Longleftrightarrow} \quad 
\left\{ \, 
\begin{array}{l} 
\un{\bu} = \un{\bv}, \\[2mm] 
\ov{\bu} = \ov{\bv}. 
\end{array}
\right. 
\end{equation*} 
  
The inclusion relation ``$\subseteq$'' in $\mbKR$ continues the inclusion relation 
in $\mbIR$ that considers intervals as sets. So, we have: 
\begin{equation}
\label{subseteq}
\bu\subseteq\bv \quad 
\stackrel{\mathrm{def}}{\Longleftrightarrow} \quad 
\left\{ \, 
\begin{array}{l} 
\un{\bu} \geq \un{\bv}, \\[2mm] 
\ov{\bu} \leq \ov{\bv}. 
\end{array} 
\right. 
\end{equation} 
In particular, $[6, 3]\subseteq[4, 5]$. 
  
The operations of taking the least upper bound (supremum) and greatest lower bound 
(infimum) with respect to inclusion are introduced for families of intervals bounded 
from above and from below respectively, using the infimum and supremum in $\mbR$: 
\begin{align*} 
\bigvee_{i\in I}\bu_i  &\  := \  \sup_{i\in I}{}_{\subseteq}\bu_i \ 
  := \;\Bigl[\;\inf_{i\in I}\un{\bu}_i,\; \sup_{i\in I}\ov{\bu}_i\;\Bigr], \\[5mm] 
\bigwedge_{i\in I}\bu_i  &\  := \  \inf_{i\in I}{}_{\subseteq}\bu_i \ 
  := \;\Bigl[\;\sup_{i\in I}\un{\bu}_i,\; \inf_{i\in I}\ov{\bu}_i\;\Bigr]. 
\end{align*} 
  
We need the following unary operations on intervals: 
\begin{align*}
\Mid\bu & := \cu := (\un{\bu} + \ov{\bu})/2 \quad\text{--- the midpoint,} \\[3mm] 
\rad\bu & := \hu := (\ov{\bu} - \un{\bu})/2 \quad\text{--- the radius,}   \\[3mm] 
\dual\bu &:= [\ov{\bu},\un{\bu}] \quad \text{--- \begin{tabular}[t]{l}
                                                 the dualization, i.\,e., swapping \\
                                                 the endpoints of the interval, 
                                                 \end{tabular}}           \\[3mm] 
\pro\bu &:= \left\{\, 
\arraycolsep=2pt 
\begin{array}{ll} 
\bu, &\text{ if } \un{\bu}\leq\ov{\bu},  \\[4pt] 
\dual\bu, &\text{ if } \un{\bu} > \ov{\bu}, 
\end{array}
\right.  \quad\text{--- \begin{tabular}[t]{l}
                        the proper projection \\ 
                        of the interval.
                        \end{tabular}} 
\end{align*}
Notice that the dualization makes sense only in $\mbKR$. 
  
Arithmetic operations of addition, subtraction, multiplication and division are 
determined through the corresponding real operations and taking exact lower and upper 
bounds by inclusion so that 
\begin{equation*} 
\arraycolsep=2pt 
\bu *\bv \  = \  \W^{\bu}\W^{\bv}(u*v), 
  \qquad \text{ where }\; \W^{\bu} := \left\{ 
\begin{array}{ll}
\bigvee\limits_{\pro\scalebox{0.75}{$\bu$}},& \mbox{ if }\;\un{\bu}\leq\ov{\bu}, \\[5mm] 
\bigwedge\limits_{\pro\scalebox{0.75}{$\bu$}},& \mbox{ if }\;\un{\bu}\geq\ov{\bu}, 
\end{array}
\right. 
\end{equation*} 
for each $*\in\{+,-,\cdot,/\}$. Naturally, division is determined only for such intervals 
$\bv$ that $0\not\in\pro\bv$. The addition and multiplication are commutative. The addition 
is defined ``by endpoints'': 
\begin{equation}
\label{sumend}
\bu + \bv = [\un{\bu} + \un{\bv}, \ov{\bu} + \ov{\bv}].
\end{equation}
The real numbers $\lambda\in\mbR$ are identified with intervals of zero radius 
$[\lambda,\lambda]$. 
Multiplication of an~interval by the number $\lambda\in\mbR$ satisfy the following 
properties: 
\begin{equation}\label{lu} 
\lambda\bu \  = \ 
\left\{\, 
\begin{array}{ll}
[\lambda\un{\bu}, \lambda\ov{\bu}], &\text{ for }\;\lambda\geq 0, \\[4pt] 
[\lambda\ov{\bu}, \lambda\un{\bu}], &\text{ for }\;\lambda\leq 0; 
\end{array} 
\right. 
\end{equation} 
\begin{equation}\label{dualu} 
(\dual\bu)\lambda\stackrel{\eqref{lu}}{=}\dual(\bu\lambda)\ 
  = \  [\ov{\bu\lambda}, \un{\bu\lambda}].
\end{equation}
The symbol $-\bu$ means the result of multiplication $(-1)\cdot\bu$, not taking 
the opposite interval for $\bu$ with respect to the addition. 
  
The matrices and vectors whose elements are intervals are called interval matrices 
and interval vectors respectively. We denote by $\bA_{i:}$ the $i$-th row of the matrix 
$\bA$. For interval vectors and matrices, the endpoints, the relations $=$ and 
$\subseteq$, the operations $\Mid$, $\rad$, $\dual$, $\pro$, as well as addition, 
subtraction and multiplication by numbers are introduced componentwise. For example, 
$(\dual\bA)_{ij}: = \dual(\bA_{ij})$, $(\bA - \bB)_{ij} := \bA_{ij} - \bB_{ij}$, 
$(-\bA)_{ij} = -\bA_{ij}$. The multiplication rules for interval vectors and matrices 
are interval extensions of analogous rules for the non-interval case: 
\begin{equation} 
\label{IMatrProd} 
(\bA\bB)_{ij}: = \sum\limits_{k} \bA_{ik}\bB_{kj}. 
\end{equation} 
  
Also, we need the property 
\begin{equation}
\label{dualAx} 
(\dual\bA)x = \dual (\bA x)
  \quad\text{ for }\bA\in\mbKR^{m\times n},\ x\in\mbR^n,
\end{equation} 
which can be easily derived from the definition of interval matrix-vector product 
(a particular case of \eqref{IMatrProd}) with the use of \eqref{sumend} and 
\eqref{dualu}.

  
\section{Results} 
  
  
\subsection{Quantifier-free descriptions in interval arithmetics}

First of all, we are going to develop quantifier-free descriptions in interval 
arithmetics for interval-quantifier linear systems and their solutions. We need 
the following notation: 
\begin{itemize} 
\item[] 
$\Qi$ will denote a quantifier prefix 
obtained from $\QAb$ by deleting all those elementary prefixes that are not related 
to the $i$-th row of the system; 
\item[] 
$\QAEAb$ will denote a quantifier prefix of the form $\QAb$ satisfying the additional 
condition: for each $i\in\{1,\ldots,m\} $ in $\Qi$, the universal quantifiers (if any) 
precede the existential quantifiers (if any); 
\item[] 
$Q^{AE}(\bA,\bb,\calA,\beta)$ will denote a quantifier prefix of the form $\QAb$ 
with the additional condition: all the universal quantifiers (if any) precede all 
the existential quantifiers (if there are such quantifiers); 
\item[] 
$\bA^\forall,\bA^\exists\in\mbIR^{m\times n}$, $\bb^\forall,\bb^\exists\in\mbIR^m$, 
$\Ac\in\mbKR^{m\times n}$, $\bc\in\mbKR^m$ will denote interval vectors and matrices 
defined by the rules 
\begin{equation}\label{AbAE}
\begin{array}{r@{\ }l@{\hspace{20mm}}r@{\ }l}
\bA^{\forall}_{\!ij}&:=\begin{cases}\bA_{ij},& \text{ if }\calA_{ij}=\forall,\\
                                   \ 0    ,& \text{ if }\calA_{ij}=\exists,
                     \end{cases}
&
\bA^{\exists}_{\!ij}&:=\begin{cases}\bA_{ij},& \text{ if }\calA_{ij}=\exists,\\
                                  \  0    ,& \text{ if }\calA_{ij}=\forall,
                     \end{cases}
\\[9mm]
\bb^{\forall}_{\!i}&:=\begin{cases}\ \bb_{i},& \text{ if }\beta_{i}=\forall,\\
                                   \ 0  ,& \text{ if }\beta_{i}=\exists,
                     \end{cases}
&
\bb^{\exists}_{\!i}&:=\begin{cases}\ \bb_{i},& \text{ if }\beta_{i}=\exists,\\
                                   \ 0  ,& \text{ if }\beta_{i}=\forall,
                     \end{cases}
\end{array}
\end{equation}
\vspace{2mm}
\begin{equation}
\label{Cd} 
\Ac_{ij}:=\begin{cases}\bA_{ij},& \text{ if }\calA_{ij}=\forall,\\ 
                      \dual\bA_{ij} ,& \text{ if }\calA_{ij}=\exists, 
                     \end{cases} 
\hspace{12mm} 
\bc_{i}:=\begin{cases}\dual\bb_{i},& \text{ if }\beta_{i}=\forall,\\ 
                      \bb_i  ,& \text{ if }\beta_{i}=\exists. 
                     \end{cases} 
\end{equation} 
\end{itemize} 
  
\bigskip 
The Gothic letter ``$\mathfrak{c}$'' as the superscript of $\mbf{A}$ and $\mbf{b}$ 
in formula \eqref{Cd} means ``characteristic''. Overall, the matrix $\Ac$ and vector 
$\bc$ will be called \emph{characteristic matrix} and \emph{characteristic vector} 
that correspond to the distribution of interval uncertainty types (A-type or E-type) 
described by the quantifier matrix $\calA$ and vector $\beta$ in the system 
under study \cite{SharyRC02,SharyBook}. 
  
We should write out the property of interval-quantifier linear systems, which we will 
repeatedly apply in the sequel: each elementary quantifier prefix from $\QAb$ can be 
carried to the row of the system in which the parameter of this prefix is present. 
This means 
\begin{equation}
\label{prop}
\QAb\ (Ax\,\sigma\,b)\qquad \Longleftrightarrow \qquad 
\conj{i\in\{1,\ldots,m\}} \Qi\ (A_{i:}x\,\sigma_i\,b_i).
\end{equation} 
The substantiation for this property is that the system of relations $(Ax\,\sigma\,b)$ 
is, in terms of logic, the conjunction of the relations, that is, 
\begin{equation*} 
\conj{i} \bigl(A_{i:} x \,\sigma_i\,b_{i}\bigr). 
\end{equation*} 
Additionally, for the conjunction, there hold equivalences  
\begin{align*}
\forall t\in{\cal S}\ \bigl(P_1(t)\ \&\ P_2\bigr)\quad &\Longleftrightarrow\quad 
\bigl(\forall t\in{\cal S}\  P_1(t)\bigr)\ \&\ P_2, \\[2mm] 
\exists t\in{\cal S}\ \bigl(P_1(t)\ \&\ P_2\bigr)\quad &\Longleftrightarrow\quad 
\bigl(\exists t\in{\cal S}\  P_1(t)\bigr)\ \&\ P_2, 
\end{align*}
where \begin{tabular}[t]{l} 
      $\cal S$ is the set of values of the variable $t$, \\[1pt] 
      $P_1$, $P_2$ are formulas, and $P_2$ does not depend on $t$. 
      \end{tabular}
  
\smallskip 
In view of \eqref{prop}, it is obvious that 
\begin{equation}\label{aeAE}
\QAEAb\ (Ax\,\sigma\,b)\quad \Longleftrightarrow\quad 
Q^{AE}(\bA,\bb,\calA,\beta)\ (Ax\,\sigma\,b),
\end{equation} 
i.\,e., the vector $x$ is a solution to the system $\QAEAb\ (Ax\,\sigma\,b)$ if and 
only if it is a~solution to the system $Q^{AE}(\bA,\bb,\calA,\beta)\ (Ax\,\sigma\,b)$. 
Thus, although the class of systems of the form $\QAEAb\,(Ax\,\sigma\,b)$ is wider 
than the class of systems of the form $Q^{AE}(\bA,\bb,\calA,\beta)\;(Ax\,\sigma\,b)$, 
the statements proved for the solutions to the system $Q^{AE}(\bA,\bb,\calA,\beta)\; 
(Ax\,\sigma\,b)$ are trivially generalized into statements for the solutions 
to the system $\QAEAb\; (Ax\,\sigma\,b)$. 
  
Now let us turn to the interval-quantifier systems of linear equations. Quantifier-free 
descrip\-tions for the widest subset of such systems have been obtained by S.P.\,Shary. 
In \cite {Shary96, Shary99}, he first proved that 
\begin{align} 
\label{Sh9699} 
Q^{AE}(\bA,\bb,\calA,\beta)\ (Ax=b)\ \ \Longleftrightarrow\ \
\bA^{\forall}x-\bb^{\forall} \subseteq \bb^\exists-\bA^{\exists}x 
\ \ \Longleftrightarrow\ \
\Ac x\subseteq\bc.
\end{align} 
(see also \cite{SharyRC02}). Equivalence \eqref{aeAE} allows us to make the following 
generalization of \eqref{Sh9699}. 
  
\begin{theorem}
\label{SharyT=} 
\begin{align}\label{Axeqb}
\QAEAb\ (Ax=b)\ \ \Longleftrightarrow\ \
\bA^{\forall}x-\bb^{\forall} \subseteq \bb^\exists-\bA^{\exists}x 
\ \ \Longleftrightarrow\ \
\Ac x\subseteq\bc.
\end{align}
\end{theorem}
  
\bigskip 
Theorem~\ref{SharyT=} for the interval-quantifier system of equations $\QAEAb\; 
(Ax = b)$ gives equivalent quantifier-free inclusion systems in $\mbIR$ 
\begin{equation*} 
\bA^{\forall} x - \bb^{\forall}\; \subseteq\; \bb^\exists - \bA^{\exists} x 
\end{equation*} 
and in $\mbKR$ 
\begin{equation*} 
\Ac x \subseteq\bc. 
\end{equation*} 
  
\bigskip\noindent
\textbf{Definition 2}\, 
Let us agree to refer to interval quantifier systems of relations, in which the vector 
of relations $\sigma$ consists of the same components, as \textit{relationally 
homogeneous systems}. 
  
\bigskip 
The results of Theorem~\ref{SharyT=} is intended for systems of equations, and our 
immediate goal is to obtain similar results for relationally homogeneous systems 
of inequalities. 
  
\begin{theorem}
\label{SharayaTineq}
\begin{align}
\QAb\ (Ax\geq b)\ \
&\Longleftrightarrow\ \
\un{\bA^{\forall}x} + \ov{\bA^{\exists}x} \geq \ov{\bb}^\forall + \un{\bb}^{\exists} 
\quad \Longleftrightarrow\quad \un{\Ac x} \geq \un{\bc},\label{Axgeqb} \\[3mm] 
\QAb\ (Ax\leq b)\ \
&\Longleftrightarrow\ \
\ov{\bA^{\forall}x} + \un{\bA^{\exists}x} \leq \un{\bb}^\forall + \ov{\bb}^{\exists} 
\quad\Longleftrightarrow\quad\ov{\Ac x} \leq \ov{\bc}. \label{Axleqb} 
\end{align}
\end{theorem}
  
\bigskip 
P r o o f. \ 
We carry out the proof only for the chain of equivalences \eqref{Axgeqb}. 
For \eqref{Axleqb}, the substantiation is completely similar. 
  
1) From \eqref{prop}, it follows that 
\begin{equation}\label{sharayaet21}
\QAb\ (Ax\geq b) \quad \Longleftrightarrow\quad 
\conj{i\in\{1,\ldots,m\}} \Qi\ (A_{i:}x\geq b_i).
\end{equation}
  
2) Using the fact that 
\begin{equation*} 
A_{i:}x\geq b_i \quad\Longleftrightarrow\quad \sum_{j=1}^n A_{ij}x_j+(-b_i)\geq 0
\end{equation*} 
and that, for any continuous functions $h:\mbR^2\to\mbR$, $g:\mbR\to\mbR$ and 
an interval $\bu\in\mbIR$, there holds 
\begin{align*}
(\forall u\in\bu)\ (h(u,x) + g(x)\geq 0)\quad &\Longleftrightarrow\quad 
\min_{u\in\scalebox{.75}{$\bu$}}\,h(u,x) + g(x)\geq 0,               \\[2mm] 
(\exists u\in\bu)\ (h(u,x) + g(x)\geq 0)\quad &\Longleftrightarrow\quad 
\max_{u\in\scalebox{.75}{$\bu$}}\,h(u,x) + g(x)\geq 0,
\end{align*}
enables us to get a quantifier-free description for \ $\Qi\ (A_{i:}x\geq b_i)$: 
\begin{equation} 
\label{sharayaet22} 
\Qi\ (A_{i:}x\geq b_i)\quad \Longleftrightarrow \quad 
  \sum_{j=1}^n \;\underset{A_{ij}\in\scalebox{.75}{$\bA$}_{ij}} 
  {\ext^{\calA_{ij}}}(A_{ij}x_j) + \underset{b_i\in\scalebox{.75}{$\bb$}_i}
  {\ext^{\beta_i}}(-b_i)\,\geq\, 0, 
\end{equation}
where ``extr'' means conditional extremum, such that 
\begin{equation*} 
\ext^{\forall} = \min, \qquad   \ext^{\exists} = \max.
\end{equation*} 
  
3) Thanks to the equalities 
\begin{equation*} 
\min_{u\in\scalebox{.75}{$\bu$}} (ux)=\underline{\bu x},
\qquad
\max_{u\in\scalebox{.75}{$\bu$}} (ux)=\overline{\bu x},
\qquad
\min_{u\in\scalebox{.75}{$\bu$}} (u)=\underline{u},
\qquad
\max_{u\in\scalebox{.75}{$\bu$}} (u)=\overline{u},
\end{equation*} 
which are valid for any interval $\bu\in\mbIR$, and taking into account \eqref{sumend}, 
the sum of the extrema in \eqref{sharayaet22} can be expressed in terms of the matrices 
$\bA^{\forall}$, $\bA^{\exists}$ and the vectors $\bb^{\forall}$, $\bb^{\exists}$ from 
\eqref{AbAE}: 
\begin{equation}\label{sharayaet23}
\sum_{j=1}^n \,\underset{A_{ij}\in\scalebox{.75}{$\bA$}_{ij}}{\ext^{\calA_{ij}}}(A_{ij}x_j) 
+\underset{b_i\in\scalebox{.75}{$\bb$}_i}{\ext^{\beta_i}}(-b_i)\geq 0
\quad \Longleftrightarrow\quad 
\underline{\bA_{i:}^{\forall}x}+\overline{\bA_{i:}^{\exists}x} \geq 
\ov{\bb}_i^\forall+\un{\bb}_i^{\exists}.
\end{equation}
  
4) From \eqref{sharayaet21}--\eqref{sharayaet23}, it follows that 
\begin{equation*} 
\QAb\ (Ax\geq b) \quad \Longleftrightarrow\quad 
  \un{\bA^{\forall}x} + \ov{\bA^{\exists}x} \geq \ov{\bb}^\forall + \un{\bb}^{\exists}. 
\end{equation*} 
  
5) Let us prove the second equivalence in the chain \eqref{Axgeqb}. For the matrix $\Ac$ 
in \eqref{Cd}, we have 
\begin{equation}\label{endsCx}
[\underline{\Ac x},\overline{\Ac x}] = \Ac x
\stackrel{\scalebox{.7}{\begin{tabular}{c}definitions \\ 
 of $\Ac,\bA^\forall,\bA^\exists$\end{tabular}}}{=}\bA^\forall x+(\dual\bA^\exists)x 
\stackrel{\scalebox{.7}{\begin{tabular}{c}properties \\ 
 \eqref{dualAx} and \eqref{sumend}\end{tabular}}}{=}
[\underline{\bA^\forall x}+\overline{\bA^\exists x},\,
  \overline{\bA^\forall x}+\underline{\bA^\exists x}],
\end{equation}
and therefore 
$\underline{\Ac x}=\underline{\bA^\forall x}+\overline{\bA^\exists x}$. The definitions 
of the vectors $\bc$, $\bb^\forall$, and $\bb^\exists$ give 
\begin{equation}\label{endsd}
[\un{\bc},\ov{\bc}]= \bc = \dual(\bb^\forall)+\bb^\exists=
[\ov{\bb}^\forall + \un{\bb}^\exists, \un{\bb}^\forall + \ov{\bb}^\exists], 
\end{equation}
hence $\un{\bc} = \ov{\bb}^\forall + \un{\bb}^\exists$. Overall, we get 
\begin{equation*} 
\underline{\bA^{\forall}x} + \overline{\bA^{\exists}x} \geq 
  \ov{\bb}^\forall+\un{\bb}^\exists \quad  \Longleftrightarrow \quad 
  \un{\Ac x} \geq \un{\bc}. 
\end{equation*} 
  
The proof of Theorem~\ref{SharayaTineq} is complete. 
  
\bigskip 
In the interval arithmetics $\mbIR$ and $\mbKR$, the relations ``$\geq$'' and 
``$\leq$'' are applicable, and they are continuations of the same relations over 
$\mbR$. For vectors, ``$\geq$'' and ``$\leq$'' are introduced componentwise. This 
allows us to formally refer to the records with $\bA^\forall$, $\bA^\exists$, 
$\bb^\forall$, $\bb^\exists$ in \eqref{Axgeqb} and \eqref{Axleqb} as inequalities 
in classical interval arithmetic, while the records with $\Ac$ and $\bc$ will be 
called inequalities in the Kaucher arithmetic. Still, in practice it is more convenient 
to understand all inequalities from \eqref{Axgeqb} and \eqref{Axleqb} as componentwise 
inequalities in $\mbR^m$. 
  
From \eqref{prop} and Theorems~\ref{SharayaTineq}, the following remarkable fact follows: 
\emph{the solution sets of interval-quantifier systems of linear inequalities with arbitrary 
$\sigma\in\{\geq, \leq \}^m$ does not depend on the order of the elementary quantifier 
prefixes, that is, all interval-quantifier systems of linear inequalities with the same 
$\bA$, $\bb$, $\cal{A}$, $\beta$ and $\sigma$ have the same solution sets}. This property
essentially distinguishes interval systems of inequalities from interval systems 
of equations. 
  
We give a corollary of Theorems~\ref{SharyT=} and \ref{SharayaTineq}, which establishes 
the relation between AE-solution sets of interval systems of linear equations and 
quantifier solution sets of interval relationally homogeneous systems of linear 
inequalities. 
  
\begin{corollary}
\label{cor1}
\begin{equation*} 
Q^{AE}(\bA,\bb,\calA,\beta)\,(Ax=b) 
\ \Longleftrightarrow\ 
\QAEAb\,(Ax=b)
\ \Longleftrightarrow\ 
\begin{cases}\QAb(Ax\geq b),\\ \QAb(Ax\leq b).\end{cases}
\end{equation*} 
\end{corollary}
  
The proof is given by the following chain of equivalences: 
\begin{multline*}
\begin{cases}\QAb\ (Ax\geq b)\\ \QAb\ (Ax\leq b)\end{cases}
\quad \stackrel{\text{Theorem~\ref{SharayaTineq}}}{\Longleftrightarrow}\quad 
\begin{cases} 
\un{\Ac x}\geq\un{\bc}\\[3pt] 
\ov{\Ac x}\leq\ov{\bc} \end{cases} 
\quad \stackrel{\text{definition of}\;\subseteq}{\Longleftrightarrow}\quad 
\Ac x\subseteq\bc 
\\[4mm] 
\stackrel{\text{Theorem~\ref{SharyT=}}}{\Longleftrightarrow}\ \ \QAEAb\ (Ax=b)
\ \ \stackrel{\eqref{aeAE}}{\Longleftrightarrow}\ \ 
Q^{AE}(\bA,\bb,\calA,\beta)\ (Ax=b).
\end{multline*}
  
\bigskip 
Theorems~\ref{SharyT=}--\ref{SharayaTineq} give quantifier-free descriptions 
for relationally homogeneous systems. Next, we turn to the consideration of systems 
with an arbitrary relationship vector $\sigma$. 
  
\bigskip\noindent
\textbf{Definition 3}\, 
We denote by $\Qsigma$ a quantifier prefix of the form $\QAb$ satisfying the following 
condition: if $\sigma_i$ is ``$=$'', then the universal quantifiers (if any) precede 
the existential quantifiers (if any) in $Q^{\sigma}_{i:} (\bA, \bb, \calA, \beta)$. 
  
\bigskip\noindent
\textbf{Definition 4}\, 
The \emph{class $Q^\sigma$} within the set of all interval-quantifier systems of linear 
relations is a subset consisting of all systems of the form $\Qsigma\ (Ax\,\sigma\,b)$. 
  
\bigskip 
The following theorem gives a quantifier-free description of the class $Q^\sigma$ 
in the interval arithmetics $\mbKR$ and $\mbIR$, with the use of componentwise 
inequalities from $\mboR^m$, where $\mboR$ denotes the extended real axis, i.\,e., 
$\mboR = \mbR\cup\{ -\infty, \infty \}$. 
  
\begin{theorem}
\label{SharayaTsigma}
\begin{equation}
\label{Cxlgd} 
\Qsigma\ (Ax\,\sigma\,b)\ \ \Longleftrightarrow\ \ 
\begin{cases}
\un{\Ac x}\geq\un{\bc} + u, \\[4pt]
\ov{\Ac x}\leq\ov{\bc} + v, \\
\end{cases}
\Longleftrightarrow\ \
\begin{cases}
\un{\bA^\forall x} + \ov{\bA^\exists x}\geq\ov{\bb}^\forall + \un{\bb}^\exists + u, \\[4pt] 
\ov{\bA^\forall x} + \un{\bA^\exists x}\leq\un{\bb}^\forall + \ov{\bb}^\exists + v, \\ 
\end{cases}
\end{equation}
where $\Ac$ and $\bc$ are from \eqref{Cd}, $\bA^\forall$, $\bA^\exists$, $\bb^\forall$, 
$\bb^\exists$ are from \eqref{AbAE}, while the vectors $u,v\in\mboR^m$ are defined as 
\begin{equation*} 
u_{i} := 
\left\{ \ 
\begin{array}{cl} 
0, &\text{if $\sigma_i$ is ``$=$'' or ``$\geq$''}, \\[4pt] 
-\infty,   &\text{if ``$\sigma_i$'' is ``$\leq$''},
\end{array}
\right. 
\qquad \ 
v_{i} := 
\left\{ \ 
\begin{array}{cl} 
0, &\text{if $\sigma_i$ is ``$=$'' or ``$\leq$''}, \\[4pt] 
\infty,   &\text{if $\sigma_i$ is ``$\geq$''}. 
\end{array} 
\right. 
\end{equation*} 
\end{theorem}
  
P r o o f \ (step by step). 
  
1) Due to the fact that each interval parameter (the element of the matrix $A$ or 
the vector $b$) enters only one row of the system $Ax\,\sigma\,b$, we have \eqref{prop} 
and, in particular, 
\begin{equation}
\label{rows} 
\Qsigma\ (Ax\,\sigma\,b)\quad \Longleftrightarrow\quad 
\conj{i\in\{1,\ldots,m\}} Q^{\sigma}_{i:}(\bA,\bb,\calA,\beta)\ (A_{i:}x\,\sigma_i\,b_i). 
\end{equation}
  
2) We eliminate quantifier prefixes in the predicate $Q^{\sigma}_{i:}(\bA,\bb,\calA,\beta)\ 
(A_{i:}x\,\sigma_i\,b_i)$ using Theorems~\ref{SharyT=} and \ref{SharayaTineq}, based on 
the specific values of $\sigma_i$: 
\begin{align*} 
Q^{\sigma}_{i:}(\bA,\bb,\calA,\beta)\ (A_{i:}x = b_i)\
&\quad \overset{\eqref{Axeqb}}{\Longleftrightarrow} \quad 
  (\Ac x)_i\subseteq\bc_i \quad  \,\Longleftrightarrow \quad 
  \bigl(\,(\un{\Ac x})_i\geq\un{\bc_i}\,\bigr)\ 
  \&\  \bigl(\,(\ov{\Ac x})_i\leq\ov{\bc_i}\,\bigr), 
  \\[5pt]
Q^{\sigma}_{i:}(\bA,\bb,\calA,\beta)\ (A_{i:}x\geq b_i)\
& \quad \overset{\eqref{Axgeqb}}{\Longleftrightarrow} \quad 
  (\un{\Ac x})_i\geq\un{\bc_i} \quad \,\Longleftrightarrow\quad 
  \bigl(\,(\un{\Ac x})_i\geq\un{\bc_i}\,\bigr)\  \&\  
  \bigl(\,(\ov{\Ac x})_i\leq\infty\,\bigr), 
  \\[5pt] 
Q^{\sigma}_{i:}(\bA,\bb,\calA,\beta)\ (A_{i:}x\leq b_i)\
&\quad \overset{\eqref{Axleqb}}{\Longleftrightarrow} \quad 
  (\ov{\Ac x})_i\leq\ov{\bc_i}\, \quad \Longleftrightarrow\quad 
  \bigl(\,(\un{\Ac x})_i\geq-\infty\,\bigr)\ \&\
  \bigl(\,(\ov{\Ac x})_i\leq\ov{\bc_i}\,\bigr). 
\end{align*}
  
3) Introducing the vectors $u$ and $v$, we pass to the matrix-vector inequalities 
\begin{equation*} 
\Qsigma\ (Ax\,\sigma\,b) 
\quad \Longleftrightarrow\quad 
\begin{cases} 
\un{\Ac x}\geq\un{\bc} + u, \\[4pt] 
\ov{\Ac x}\leq\ov{\bc} + v. \\
\end{cases}
\end{equation*} 
  
4) The equivalence 
\begin{equation*} 
\begin{cases}
\un{\Ac x}\geq\un{\bc} + u, \\[4pt] 
\ov{\Ac x}\leq\ov{\bc} + v, \\
\end{cases}
\quad \Longleftrightarrow\quad 
\begin{cases}
\un{\bA^\forall x} + \ov{\bA^\exists x} \geq\ov{\bb}^\forall + \un{\bb}^\exists + u, \\[4pt] 
\ov{\bA^\forall x} + \un{\bA^\exists x} \leq\un{\bb}^\forall + \ov{\bb}^\exists + v, \\ 
\end{cases} 
\end{equation*} 
is obvious in view of \eqref{endsCx} and \eqref{endsd}. 
The proof of Theorem~\ref{SharayaTsigma} is complete. 
\medskip
  
Convenient quantifier-free representations for the class $Q^\sigma$ can be obtained 
from Theorem~3, if we introduce the sets of intervals $\mbKoR = \{[\un{z}, \ov{z}]\mid 
\un{z}, \ov{z}\in\mboR \}$ and $\mbIoR = \{[\un{z}, \ov{z}] \mid \un{z}, \ov{z} \in \mboR, 
\  \un{z}\leq\ov{z} \}$ and continue relation ``$\subseteq$'' according to rule 
\eqref{subseteq}. Then 
\begin{equation}\label{Cxinds}
\Qsigma\ (Ax\,\sigma\,b)\quad \Longleftrightarrow\quad
  \Ac x\subseteq \bc + \bw \quad \Longleftrightarrow\quad 
  \bA^\forall x - \bb^\forall\subseteq\bb^\exists - \bA^\exists x + \bw, 
\end{equation} 
where $\Ac$ and $\bc$ from \eqref{Cd}, $\bA^\forall$, $\bA^\exists$, $\bb^\forall$, 
$\bb^\exists$ from \eqref{AbAE}, and the interval vector $\bw\in\mbIoR^m $ is such that 
\begin{equation*} 
\bw_i:=\left\{\begin{array}{c@{\hspace{3mm}}c}
0,             & {\text{ if $\sigma_i$ is ``$=$''}},   \\[3pt]
{[0,\infty]},  & {\text{ if $\sigma_i$ is ``$\geq$''}},\\[3pt] 
{[-\infty,0]}, & {\text{ if $\sigma_i$ is ``$\leq$''}}.
\end{array}\right.
\end{equation*} 
  
The inclusion 
\begin{equation*} 
\Ac x\;\subseteq\;\bc + \bw 
\end{equation*} 
provides a quantifier-free description of the solution set to the quantifier interval 
linear system $\Qsigma \ (Ax\,\sigma\,b)$ in any interval arithmetic that extends 
the Kaucher arithmetic to the set $\mbKoR$. An example of such an extension is given 
in \cite{Kaucher73}. We agree to denote the arithmetic extension, as well as its basic 
set, through $\mbKoR$. 
  
Similarly, the inclusion 
\begin{equation*} 
\bA^\forall x - \bb^\forall \;\subseteq\; \bb^\exists - \bA^\exists x + \bw
\end{equation*} 
provides a quantifier-free description of the solution set to the system $\Qsigma \  
(Ax\,\sigma\,b)$ in interval arithmetic that extends $\mbIR$ to the set $\mbIoR$. 
Examples of the extension of the classical interval arithmetic to a set of intervals 
with infinite endpoints are described in \cite{Markov92}. Let us agree to refer to 
any such extension as arithmetic $\mbIoR$. Thus, the relation \eqref{Cxinds} gives 
quantifier-free descriptions of the solution sets to quantifier interval linear 
systems of class $Q^\sigma$ in the interval arithmetics $\mbKoR$ and $\mbIoR$. 
  
Comparing the quantifier-free descriptions obtained for the solution sets to quantifier 
interval linear systems, we can see that, 
\begin{description} 
\item 
on the one hand, the quantifier-free description in $\mbKR$ ($\mbKoR$) is much more 
remote from the initial data $\bA$, $\bb$, $\cal A$ and $\beta$ due to multilevel 
notation, and, 
\item 
on the other hand, the description in $\mbKR$ ($\mbKoR$) is more concise and 
convenient for analysis than a similar description in $\mbIR$ ($\mbIoR$). 
\end{description}

  
\subsection{Quantifier-free descriptions in real arithmetic}

In this section, we derive quantifier-free descriptions of the quantifier solution 
sets to interval linear systems in the real arithmetic $\mbR$. To do that, we will 
need \emph{Hadamard product} of matrices (entrywise product), denoted by the symbol 
``$\circ$'' (see e.\,g. \cite{HornJohnson}). Hadamard product is defined for 
two matrices of the same dimensions and produces another matrix in which the $ij$-th 
element is the product of the $ij$-th elements of the original matrices: 
\begin{equation*} 
(A\circ B)_{ij} = A_{ij}B_{ij}. 
\end{equation*} 
Also, notice that the operation of taking the modulus of a vector is understood 
componentwise. If, for instance, $x\in\mbR^n$, then $|x|$ is a nonnegative vector 
with the components $|x|_i = |x_i|$. 
  
\begin{theorem}
\label{SharayaTR}
\begin{align}
\QAEAb\ (Ax=b) &\quad \Longleftrightarrow\quad 
|\cA x-\cb|\leq(\cAs\!\circ\!\hA)|x|+\betas\!\circ\,\hb,\label{TR=}\\
\QAb\ (Ax\geq b) &\quad \Longleftrightarrow\quad 
\cb-\cA x\,\leq(\cAs\!\circ\!\hA)|x|+\betas\!\circ\,\hb,\label{TRgeq}\\
\QAb\ (Ax\leq b) &\quad \Longleftrightarrow\ \quad
\cA x-\cb\,\leq(\cAs\!\circ\!\hA)|x|+\betas\!\circ\,\hb,\label{TRleq}\\
\Qsigma\ (Ax\,\sigma\,b) &\quad \Longleftrightarrow\quad 
\abs^{\sigma}(\cA x-\cb)\leq(\cAs\!\circ\!\hA)|x|+\betas\!\circ\,\hb,\label{TRs}
\end{align}
where 
\begin{equation}\label{Asbs}
\cAs_{ij}=\begin{cases}\ \ 1, & \text{ if }\calA_{ij}=\exists,\\
                        -1, & \text{ if }\calA_{ij}=\forall, \end{cases}
\qquad
\betas_i=\begin{cases} \ \ 1, & \text{ if }\beta_i=\exists,\\
                        -1, & \text{ if }\beta_i=\forall, \end{cases}
\end{equation}
\begin{equation*} 
\abs_i\!^{\sigma}(y) = 
\left\{ \ 
\begin{array}{cc}
|y_i|, & \text{ if $\sigma_i$ is  ``$=$''},   \\[3pt]
 -y_i, & \text{ if $\sigma_i$ is ``$\geq$''}, \\[3pt]
 y_i,  & \text{ if $\sigma_i$ is ``$\leq$''}. 
\end{array} 
\right. 
\end{equation*} 
\end{theorem}
  
P r o o f. 
  
1) The equivalence \eqref{TR=} was proposed and proved by Jiri Rohn at the international 
conference INTERVAL'96 (September--October of 1996, W\"{u}rzburg, Germany), in a private 
talk with Sergey Shary and Anatoly Lakeyev. Later, its reformulation with the use 
of Hadamard product was proposed by Anatoly Lakeyev in the work \cite{Lakeyev03}. 
Below, we present our own proof.

In view of Theorem~\ref{SharyT=}, 
\begin{equation*} 
\QAEAb (Ax=b)\quad \Longleftrightarrow\quad \Ac x\subseteq\bc.
\end{equation*} 
Then, using the properties of Kaucher arithmetic 
\begin{equation*} 
(\forall\bu,\bv\in\mbKR^m)\ \
\bigl(\bu\subseteq\bv\ \Longleftrightarrow\ |\cu-\cv|\leq\hv-\hu\bigr),
\end{equation*} 
\begin{equation}
\label{midradCx}
\Mid(\Ac x) = \check{\Ac} x,\qquad
\rad(\Ac x) = \hat{\Ac}\,|x|,
\end{equation}
and we get 
\begin{equation*} 
\Ac x\subseteq\bc \quad \Longleftrightarrow\quad 
  |\check{\Ac} x- \check{\bc}|\;\leq\;\hat{\bc} - \hat{\Ac}\,|x|. 
\end{equation*} 
From the definitions of \eqref{Cd} and \eqref{Asbs} for $\Ac$, $\bc$, $\cAs$, 
and $\betas$, we have 
\begin{equation}\label{ch}
\check{\Ac} = \cA, 
  \qquad \hat{\Ac} = -\cAs\circ\hA, 
  \qquad \check{\bc}=\cb, 
  \qquad \hat{\bc} = \betas\circ\hb.
\end{equation}

2) Let us prove the equivalence \eqref{TRgeq}. According to Theorem~\ref{SharayaTineq} 
\begin{equation*} 
\QAb (Ax\geq b)\quad \Longleftrightarrow\quad \underline{\Ac x}\;\geq\;\un{\bc}.
\end{equation*} 
Drawing on the obvious property of the Kaucher arithmetic 
\begin{equation*} 
(\forall\bu,\bv\in\mbKR^m)\ \
\bigl(\un{u}\geq\un{v}\ \Longleftrightarrow\  \cv-\cu\leq\hv-\hu\bigr),
\end{equation*} 
which allows us to replace the inequality between the endpoints by the inequality 
between centers and radii, and then involving \eqref{midradCx}, we get 
\begin{equation*} 
\underline{\Ac x}\,\geq\,\un{\bc}
  \quad \Longleftrightarrow\quad 
  \check{\bc} - \check{\Ac} x \,\leq\, \hat{\bc} - \hat{\Ac}\,|x|. 
\end{equation*} 
Finally, we use \eqref{ch}. 
  
3) The equivalence \eqref{TRleq} is proved similarly to \eqref{TRgeq}. 
  
4) It remains to substantiate the equivalence \eqref{TRs}. Just as in the item~1 
of the proof of Theorem~\ref{SharayaTsigma}, we have \eqref{rows}, i.\,e., the problem 
splits in rows. We apply, to each row, one of the equivalences \eqref{TR=}, \eqref{TRgeq} 
or \eqref{TRleq}, depending on the corresponding binary relation, and convolve 
the resulting system of inequalities using the operation $\abs^\sigma$. 
  
The proof of Theorem~\ref{SharayaTR} is complete. 
\medskip
  
From the equivalences \eqref{TR=}--\eqref{TRleq}, one more proof of Corollary~\ref{cor1} 
becomes obvious. In addition, it is not difficult to establish the following connection 
between relationally homogeneous systems of inequalities of the opposite signs. 
  
\begin{corollary}
\label{cor2}
\begin{gather}
\QAb\ (Ax\geq b) \quad\Longleftrightarrow\quad 
  Q(-\bA,-\bb,\calA,\beta)\ (Ax\leq b),\label{geqleq1} \\[2mm] 
Q(-\bA,-\bb,\calA,\beta)\ (Ax\geq b) \quad\Longleftrightarrow\quad 
  \QAb\ (Ax\leq b).\label{geqleq2}
\end{gather}
\end{corollary}
  
P r o o f. \  Based on the properties of intervals 
\begin{equation}\label{midrad-}
\Mid(-\bu)=-\cu\qquad \text{ and }\qquad \rad(-\bu)=\hu,
\end{equation} 
we can show the validity of relation~\eqref{geqleq2}: 
\begin{align*}
Q(-\bA,-\bb,\calA,\beta)\ (Ax\geq b) &\ \stackrel{\eqref{TRgeq}}{\Longleftrightarrow} \ 
  \Mid(-\bb)-\Mid(-\bA)x\,\leq\bigl(\cAs\!\circ\,\rad(-\bA)\bigr)|x| 
  + \betas\!\circ\,\rad(-\bb),                               \\[2mm] 
&\ \stackrel{\eqref{midrad-}}{\Longleftrightarrow}\ \ 
  -\cb+\cA x\,\leq(\cAs\!\circ\!\hA)|x|+\betas\!\circ\,\hb,  \\[2mm] 
&\ \stackrel{\eqref{TRleq}}{\Longleftrightarrow}\ \  \QAb\ (Ax\leq b). 
\end{align*}
Relation~\eqref{geqleq1} is proved similarly. 
  
\medskip 
Corollary~\ref{cor2} means that, if the sign of the inequality and the signs of all 
intervals of the parameter values are reversed to the opposite, then the set of quantifier 
solutions to the interval system of linear inequalities does not change. For example, 
the sets of solutions to the systems $(\forall A\in\bA)\,(\exists b\in\bb)\,(Ax\geq b)$ 
and $(\forall A\in -\bA)\,(\exists b\in -\bb)\,(Ax\leq b)$ coincide. 
\bigskip
  
  
\subsection{Quantifier-free descriptions in $\mbKR$, $\mbIR$ and $\mbR$ \\ 
            for systems of basic types}

So far, considering the interval-quantifier linear systems $\QAb\,(Ax\,\sigma\,b)$, 
we tried to obtain results in which there are no constraints on the parameters 
$\bA$, $\bb$, $\calA$, $\beta$ and $\sigma$, and the restrictions on the order of 
the elementary quantifier prefixes in $Q$ are minimal. In this sense, the most general 
descriptions were found for the class $Q^{\sigma}$. In this section, we consider subsets 
of interval-quantifier linear systems of class $Q^{\sigma}$, which are distinguished by 
the requirement of homogeneity of $\calA$ and the homogeneity of $\beta$. Elements of 
all these subsets will be called systems of \textit{basic types}, and their solutions 
will be referred to as \emph{quantifier solutions of basic types} for interval linear 
systems of the form $\bA x\,\sigma\,\bb$. 
  
Depending on which quantifiers fill the matrix $\cal{A}$ and the vector $\beta$, all 
the interval-quantifier linear systems of the basic types are divided into 4 subsets, 
or 4 types. This subdivision is presented in the last column of Table~1. For solutions 
to systems of each of the main types, we give a proper name that continues the one used 
in \cite{FNRRZ,SharyBook} for solutions of relationally homogeneous systems of this type. 
The names of the solutions are listed in the first column of Table~1, the values of 
the elements of the matrix $\cal{A}$ and components of the vector $\beta$ are listed 
in the second and third columns, and the fourth column gives the general form for 
the interval-quantifier systems of the corresponding basic type. 
  
  
\begin{table}[hbt]
\begin{center}
\begin{tabular}{|l|c|c|c|}
\multicolumn{4}{r}{T\ a\ b\ l\ e\ \ 1}\\[1mm]
\multicolumn{4}{c}{{\bf Basic types of quantifier solutions to the interval system 
$\bA x\,\sigma\,\bb$}}\\[2mm]
\hline
\multicolumn{1}{|c|}{\raisebox{-1mm}{Name}} & \multicolumn{2}{c|}{Values of elements} 
&\multicolumn{1}{c|}{\raisebox{-1mm}{Interval-quantifier}}\\ \cline{2-3}
\multicolumn{1}{|c|}{of solution} & of the matrix $\cal{A}$ & of the vector $\beta$ 
& \multicolumn{1}{c|}{system of basic type}\\ \hline\hline
& & &\\[-3mm]
weak 
& $\exists$ & $\exists$ 
& $(\exists A\in\bA) (\exists b\in\bb)\ (Ax\,\sigma\,b)$\\[3pt] 
tolerable 
& $\forall$ & $\exists$ 
& $(\forall A\in\bA) (\exists b\in\bb)\ (Ax\,\sigma\,b)$\\[3pt] 
controllable 
& $\exists$ & $\forall$ 
& $(\forall b\in\bb) (\exists A\in\bA)\ (Ax\,\sigma\,b)$\\[3pt] 
strong 
& $\forall$ & $\forall$ 
& $(\forall A\in\bA) (\forall b\in\bb)\ (Ax\,\sigma\,b)$ \\[3pt] 
\hline
\end{tabular}
\end{center}
\end{table}
  
  
Quantifier-free descriptions in $\mbKR$, $\mbIR$ and $\mbR$ for systems of basic types 
can be obtained as corollaries of the corresponding descriptions for systems of class 
$Q^\sigma$. Let us explain that for relationally homogeneous systems using Table~2. 
  
  
\begin{table*}[tbp] 
\scalebox{.93}{
\rotatebox{90}{%
\renewcommand{\arraystretch}{2}
\begin{tabular}{||l|c||c|c|c|c|c||}
\multicolumn{7}{r}{}\\[5mm]
\multicolumn{7}{r}{T\ a\ b\ l\ e\ \  2}\\ 
\multicolumn{7}{c}{\bf Characterization of relationally homogeneous
interval-quantifier linear systems and their main solution types}\\[2mm]
\hline\hline
&  
& \multicolumn{5}{c||}{Type of solution and corresponding quantifier prefix $\QAb$}
\\ \cline{3-7} 
\multicolumn{1}{||c|}{$Ax\,\sigma\,b$}
& \rotatebox{90}{\makebox[0mm]{\begin{tabular}{@{}c@{}}\\[-13mm] Space\\[-6mm] 
of description\\[-1mm]\end{tabular}}}
&  & \multicolumn{4}{c||}{Basic types of solutions}\\ \cline{4-7}
& & Quantifier & Weak & Tolerable & Controllable & Strong\\[-3mm] 
& 
& $\Qsigma$ & $(\exists A\in\bA)\ (\exists b\in\bb)$ 
& $(\forall A\in\bA)\ (\exists b\in\bb)$ & $(\forall b\in\bb)\ (\exists A\in\bA)$
& $(\forall A\in\bA)\ (\forall b\in\bb)$\\
\hline\hline
$Ax=b$ & $\mbKR$ & $\Ac x\subseteq\bc$  & $(\dual\bA)x\subseteq\bb$ 
& $\bA x\subseteq\bb$ & $(\dual\bA)x\subseteq\dual\bb$ & $\bA x\subseteq\dual\bb$
\\ \cline{2-7} 
& $\mbIR$ & $\bA^\forall x-\bb^\forall\subseteq\bb^\exists-\bA^\exists x$  
& $0\in\bb-\bA x$ & $\bA x\subseteq\bb$ & $\bb\subseteq\bA x$ & $\bA x-\bb\subseteq0$
\\ \cline{2-7} 
& $\mbR$ & $|\cA x-\cb|\leq(\cAs\!\circ\!\hA)|x|+\betas\!\circ\!\hb$
& $|\cA x-\cb|\leq \hA|x|+\hb$ & $|\cA x-\cb|\leq -\hA|x|+\hb$ 
& $|\cA x-\cb|\leq \hA|x|-\hb$ & $|\cA x-\cb|\leq -\hA|x|-\hb$\\ \hline\hline
$Ax\geq b$ 
& $\mbKR$ 
& $\un{\Ac x}\geq\un{\bc}$  
& $\un{(\dual\bA)x}\geq\un{\bb}$ 
& $\un{\bA x}\geq\un{\bb}$ 
& $\un{(\dual\bA)x}\geq\ov{\bb}$ 
& $\un{\bA x}\geq\ov{\bb}$
\\ \cline{2-7} 
& $\mbIR$ 
& $\ov{\bA^\exists x} + \un{\bA^\forall x}\geq\un{\bb}^\exists + \ov{\bb}^\forall$ 
& $\ov{\bA x}\geq\un{\bb}$ 
& $\un{\bA x}\geq\un{\bb}$ 
& $\ov{\bA x}\geq\ov{\bb}$ 
& $\un{\bA x}\geq\ov{\bb}$ 
\\ \cline{2-7} 
& $\mbR$ & 
$\cb-\cA x\leq(\cAs\!\circ\!\hA)|x|+\betas\!\circ\!\hb$
& $\cb-\cA x\leq \hA|x|+\hb$ & $\cb-\cA x\leq -\hA|x|+\hb$ 
& $\cb-\cA x\leq \hA|x|-\hb$ & $\cb-\cA x\leq -\hA|x|-\hb$\\ \hline\hline
$Ax\leq b$ 
& $\mbKR$ 
& $\ov{\Ac x}\leq\ov{\bc}$ 
& $\ov{(\dual\bA)x}\leq\ov{\bb}$ 
& $\ov{\bA x}\leq\ov{\bb}$ 
& $\ov{(\dual\bA)x}\leq\un{\bb}$ 
& $\ov{\bA x}\leq\un{\bb}$ 
\\ \cline{2-7} 
& $\mbIR$ 
& $\un{\bA^\exists x} + \ov{\bA^\forall x}\leq\ov{\bb}^\exists + \un{\bb}^\forall$ 
& $\un{\bA x}\leq\ov{\bb}$ 
& $\ov{\bA x}\leq\ov{\bb}$ 
& $\un{\bA x}\leq\un{\bb}$ 
& $\overline{\bA x}\leq\un{\bb}$ 
\\ \cline{2-7} 
& $\mbR$ 
& $\cA x-\cb\leq(\cAs\!\circ\!\hA)|x|+\betas\!\circ\!\hb$
& $\cA x-\cb\leq \hA|x|+\hb$ & $\cA x-\cb\leq -\hA|x|+\hb$ 
& $\cA x-\cb\leq \hA|x|-\hb$ & $\cA x-\cb\leq -\hA|x|-\hb$\\ \hline\hline
\end{tabular}\renewcommand{\arraystretch}{1}}
}\end{table*} 
  
  
In the Table~2, columns~4--7, corresponding to the basic types of quantifier solutions, 
are obtained, in row-wise manner, from column~3 corresponding to quantifier solutions 
with the prefix $Q^\sigma$. It is necessary to use the definition of \eqref{Cd} 
of the matrix $\Ac$ and vector $\bc$ in the rows corresponding to the Kaucher 
arithmetic. In the rows that correspond to the classical interval arithmetic, we have 
to use definition \eqref{AbAE} of the matrices $\bA^{\forall}$, $\bA^{\exists}$ 
and the vectors $\bb^{\forall}$, $\bb^{\exists}$. Finally, the rows corresponding 
to real non-interval arithmetic, the definition \eqref{Asbs} of the matrix $\cAs$, 
vector $\betas$ and the definition of the product $\circ$ should be used. 
  
Approximately half of the descriptions of the basic types of quantifier solutions 
for interval linear systems, presented in columns 4--7 of Table 2, have been obtained 
earlier. The descriptions that were found first, obtained their own proper names. 
These are 
\begin{itemize} 
\item[] 
the Oettli-Prager characterization in $\mbR$ \cite{OettliPrager} and the Beeck 
characterization in $\mbIR$ \cite{Beeck} \\  for weak solutions of the equation 
$\bA x = \bb$, 
\item[] 
the Gerlach description in $\mbR$ for weak solutions of the inequality $\bA x\leq\bb$ 
\cite{Gerlach}. 
\end{itemize} 
The quantifier-free descriptions of the set of tolerable solutions to the equation 
$\bA x = \bb$ was obtained in $\mbR$ by J.\,Rohn \cite{Rohn85} and in $\mbIR$ by 
A.\,Neumaier \cite{Neumaier86}. The description in $\mbR$ was further investigated 
by A.V.\,Lakeyev and S.I.\,Noskov in \cite{LakeyevNoskov94}, and they also presented, 
as an evident one, a description, in $\mbIR$, for the set of controllable solutions to 
the equation $\bA x = \bb$ (see also \cite{Shary92}). The remaining descriptions for 
the basic types of quantifier solutions to the equation $\bA x = \bb$ in the interval 
arithmetics $\mbIR$ and $\mbKR$ are also known, for example, as obvious corollaries 
of the statement~\eqref{Sh9699}, proved by S.P\,Shary in \cite{Shary96,Shary99}. 
In Theorem~2.25 from the book \cite{FNRRZ}, a quantifier-free description in $\mbR$ 
for strong solutions to the interval inequality $\bA x\leq\bb$ was presented. 
  
For interval-quantifier systems of basic types in which the relationship vector 
$\sigma$ is not homogeneous, quantifier-free descriptions in $\mbKoR$ can be obtained 
from \eqref{Cxinds} and \eqref{Cd}. The descriptions in $\mbIoR$ can be derived from 
\eqref{Cxinds} and \eqref{AbAE}, and the descriptions in $\mbR$ follows from 
\eqref{TRs} and \eqref{Asbs}. Below, we give these descriptions only in $\mbIoR$ 
and $\mbR$ (in $\mbKoR$, they are less expressive and differ from the descriptions 
in $\mbIoR$ by obvious arithmetic transformations, in the same way as the descriptions 
in $\mbKR$ and $\mbIR$ differ from each other in Table~2): 
\begin{align*}
(\exists A\in\bA) (\exists b\in\bb)\ &(Ax\,\sigma\,b)
\ \ \Longleftrightarrow
&&0\in\bb-\bA x+\bw
&&\hspace{-3mm}\Longleftrightarrow \  \ 
\abs^{\sigma}(\cA x-\cb)\leq\hA|x|+\hb;
\\[3pt] 
(\forall A\in\bA) (\exists b\in\bb)\ &(Ax\,\sigma\,b)
\ \ \Longleftrightarrow
&&\bA x\subseteq\bb+\bw
&&\hspace{-3mm}\Longleftrightarrow \  \ 
\abs^{\sigma}(\cA x-\cb)\leq-\hA|x|+\hb;
\\[3pt] 
(\forall b\in\bb) (\exists A\in\bA)\ &(Ax\,\sigma\,b)
\ \ \Longleftrightarrow
&&\bb\subseteq\bA x+\bw
&&\hspace{-3mm}\Longleftrightarrow \  \ 
\abs^{\sigma}(\cA x-\cb)\leq\hA|x|-\hb;
\\[3pt] 
(\forall A\in\bA) (\forall b\in\bb)\ &(Ax\,\sigma\,b)
\ \ \Longleftrightarrow
&&\bA x-\bb\subseteq\bw
&&\hspace{-3mm}\Longleftrightarrow \  \ 
\abs^{\sigma}(\cA x-\cb)\leq-\hA|x|-\hb. 
\end{align*} 
  
  
\section{Conclusion}

The main results of the paper are presented in Theorems~\ref{SharayaTineq}--\ref{SharayaTR} 
(equivalence \eqref{TR=} was previously known) and in Corollary~\ref{cor1}. 
  
Among the statements that have no restrictions on the parameters $\bA$, $\bb$, $\calA$, 
$\beta$ and $\sigma$, those that give quantifier-free descriptions of interval-quantifier 
linear systems of class $Q^\sigma$ have the greatest generality. These are the relation 
\eqref{Cxlgd}, which provides a transition to $\mbKR$ and $\mbIR$, the relation 
\eqref{Cxinds} for the transition to $\mbKoR$ and $\mbIoR$, and equivalence \eqref{TRs} 
that allows us to go into $\mbR$. 
  
The usefulness of quantifier-free descriptions from \eqref{Cxlgd}, \eqref{Cxinds} and 
\eqref{TRs} is that they give us the possibility 
\begin{itemize} 
\item 
to study all interval-quantifier linear systems of class $Q^\sigma$ simultaneously 
and in a uniform way, and to derive results for their subclasses (in particular, for 
interval-quantifier systems of basic types) as consequences of the general result; 
\item 
to design such solution methods for problems related to interval-quantifier linear systems 
that are suitable for all systems of class $Q^\sigma$ (an example is the author's software 
package for visualization of quantifier solution sets to interval linear systems, available 
at \url{http://www.nsc.ru/interval/sharaya}). 
\end{itemize} 
  
Quantifier-free descriptions, in interval arithmetic, for various classes of 
interval-quantifier linear systems and for their solutions, both previously known 
(for example, relation  \eqref{Sh9699}) and those obtained in this paper in the form 
of relations \eqref{Axeqb}--\eqref{Axleqb}, \eqref{Cxlgd}, \eqref{Cxinds}, allow us 
\begin{itemize} 
\item 
to investigate interval-quantifier linear systems by interval methods, i.\,e., to reveal 
the properties of their solution sets, the relationships between systems with various 
conditions on the parameters $\bA$, $\bb$, $\calA$, $\beta$, $\sigma$ and the order 
of the quantifier prefixes (an example is the proof of Corollary~\ref{cor1}); 
\item 
to construct interval (that is, essentially using interval arithmetic) solution methods 
for problems in which the formulation involves interval-quantifier linear systems 
(examples of such methods for systems of equations can be found in \cite{SharyBook}, 
while for inequalities and systems of class $Q^\sigma$ constructing such methods is 
a matter of the future). 
\end{itemize}


  
\end{document}